\documentclass[a4paper,11pt]{article}
\usepackage{graphicx}
\usepackage{fancybox}
\usepackage{amssymb}
\usepackage[all]{xy}
\usepackage{color}

\def\rien{\rule{0pt}{0pt}}
\begin{document}
%
%
\title{Cross-Entropic Learning of a Machine for the Decision in a Partially Observable Universe}
%
%
\author{Fr\'ed\'eric Dambreville\\
D\'el\'egation G\'en\'erale pour l'Armement, DGA/CTA/DT/GIP\\
16 Bis, Avenue Prieur de la C\^ote d'Or\\
F 94114, France\\
Email: {\tt http://email.fredericdambreville.com}}
%
\maketitle
\begin{abstract}
In this paper, we are interested in optimal decisions in a partially observable universe.
Our approach is to directly approximate an optimal strategic tree depending on the observation.
This approximation is made by means of a parameterized probabilistic law.
A particular family of hidden Markov models, with input \emph{and} output, is considered as a model of policy.
A method for optimizing the parameters of these HMMs is proposed and applied.
This optimization is based on the cross-entropic principle for rare events simulation developed by Rubinstein.
\end{abstract}
%
\newcommand{\keywordsname}{Keywords}
\newenvironment{keywords}%
  {\small
    \list{}{\labelwidth0pt
      \leftmargin0pt \rightmargin\leftmargin
      \listparindent\parindent \itemindent0pt
      \parsep0pt
      \let\fullwidthdisplay\relax}%
    \item[\hskip\labelsep\bfseries\keywordsname:]}{\endlist}

\paragraph{Notations.}
Some specific notations are used in this document.
\begin{itemize}
\item The variables $d$, $y$, $x$ and $m$ are used for the decision, observation, world state and machine memory,
\item The time $t$ is starting from stage $1$ to the maximal stage $T$.
Variables with subscript outside this scope are synonymous to $\emptyset$.
For example, $\prod_{t=1}^T\pi(x_t|x_{t-1})$ means $\pi(x_1|\emptyset)\prod_{t=2}^T\pi(x_t|x_{t-1})$\,, \emph{ie.} a Markov chain.
A similar principle is used for the \emph{level} supscript $\lambda$ in the definition of hierarchical HMM,
\item The generic notation for a probability is $P$. However, the functions $p$, $\pi$ and $h$ denote some specific components of the probability.
$p$ is the law of the observation $y$ and state $x$ conditionally to the decision $d$.
$\pi$ is a stochastic policy, \emph{ie.} a law of the decision conditionally to the observation.
$h$ is an approximation of $\pi$ by a HMM family.
The hidden state of $h$ is defined as the machine memory $m$.
\end{itemize}
\section{Introduction}
There are different degrees of difficulty in planning and control problems.
In most problems, the planner have to start from a given state and terminate in a required final state.
There are several transition rules, which condition the sequence of decision.
For example, a robot may be required to move from room A, starting state, to room B, final state; its decision could be \emph{go forward}, \emph{turn right} or \emph{turn left}, and it cannot cross a wall; these are the conditions over the decision.
A first degree in the difficulty is to find at least one solution for the planning.
When the states are only partially known or the resulting actions are not deterministic, the difficulty is quite enhanced:
the planner has to take into account the various observations.
Now, the problem becomes much more complex, when this planning is required to be optimal or near-optimal.
For example, find the shortest trajectory which moves the robot from room A to room B.
There are again different degrees in the difficulty, depending on the problem to be deterministic or not, depending on the model of the future observations.
In the particular case of a Markovian problem with the full observation hypothesis, the dynamic programming principle\cite{IEEEisda:bellman} could be efficiently applied (Markov Decision Process theory/MDP).
This solution has been extended to the case of partial observation (Partially Observable Markov Decision Process/POMDP), but this solution is generally not practicable, owing to the huge dimension of the variables\cite{IEEEisda:sondik,IEEEisda:rocco}.
\\\\
For such reason, different methods for approximating this problem has been introduced.
For example, Reinforcement Learning methods \cite{sutton} are able to learn an evaluation table of the decision conditionnally to the known universe states and an observation short range.
In this case, the range of observation is indeed limited in time, because of an exponential grow of the table to learn.
Recent works\cite{IEEE:schmidhuber} are investigating the case of hierarchical RL, in order to go beyond this range limitation.
Whatever, these methods are generally based on an additivity hypothesis about the reward.
Another viewpoint is based on the direct learning of the policy\cite{meuleau}.
Our approach is of this kind.
It is particularly based on the Cross-Entropy optimisation algorithm developed by Rubinstein\cite{CEbook}.
This simulation method relies both on a probabilistic modelling of the policies (in this paper, these models are Bayesian Networks) and on an efficient and robust iterative algorithm for optimizing the model parameters.
More precisely, the policy will be modelled by conditional probabilistic law, \emph{i.e.} decisions depending on observations, which are involving memories; typically hidden Markov models are used.
Also are implemented a hierachical modelling of the policies by means of hierarchical hidden Markov models.
\\\\
The next section introduces some formalism and gives a quick description of the optimal planning in partially observable universes.
It is proposed a near-optimal planning method, based on the direct approximation of the optimal decision tree.
The third section introduces the family of Hierarchical Hidden Markov Models being in use for approximating the decision trees.
The fourth section describes the method for optimizing the parameters of the HHMM, in order to approximate the optimal decision tree for the POMDP problem.
The cross-entropy method is described and applied.
The fifth section gives an example of application.
A comparison with a Reinforcement Learning method, the Q-learning, is made.
The paper is then concluded.
\section{Decision in a partially observable universe}
\begin{figure}
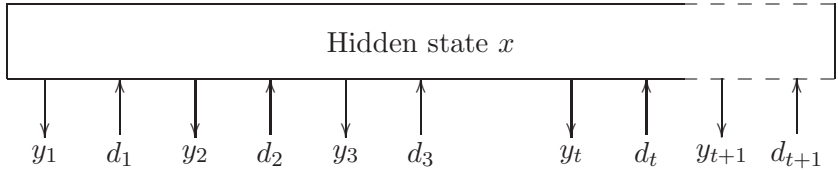

\caption{The world}
\label{ISDA:fig:1}
\begin{center}
\begin{tabular}{c}
\vspace{-25pt}\\
\vspace{-15pt}\rien\hspace{-30pt}\scalebox{1}{\xy
(-12.5,2.5)*+{\rien},
(0,0)="box00",
(90,0)="box10",
(90,-10)="box11",
(110,0)="box10f",
(110,-10)="box11f",
(0,-10)="box01",
\ar @{-}"box00";"box10",
\ar @{--}"box10";"box10f",
\ar @{-}"box10f";"box11f",
\ar @{--}"box11";"box11f",
\ar @{-}"box11";"box01",
\ar @{-}"box00";"box01",\POS "box01",
(55,-5)*+{\mbox{Hidden state $x$}},
(5,-10)="a10y",
(25,-10)="a11y",
(45,-10)="a12y",
(75,-10)="a13y",
(95,-10)="a14y",
(15,-10)="a10x",
(35,-10)="a11x",
(55,-10)="a12x",
(85,-10)="a13x",
(105,-10)="a14x",
(5,-20)*+{y_1}="b10y",
(25,-20)*+{y_2}="b11y",
(45,-20)*+{y_3}="b12y",
(75,-20)*+{y_t}="b13y",
(95,-20)*+{y_{t+1}}="b14y",
(15,-20)*+{d_1}="b10x",
(35,-20)*+{d_2}="b11x",
(55,-20)*+{d_3}="b12x",
(85,-20)*+{d_t}="b13x",
(105,-20)*+{d_{t+1}}="b14x",
\ar @{->}"a10y";"b10y",
\ar @{->}"a11y";"b11y",
\ar @{->}"a12y";"b12y",
\ar @{->}"a13y";"b13y",
\ar @{->}"a14y";"b14y",
\ar @{->}"b10x";"a10x",
\ar @{->}"b11x";"a11x",
\ar @{->}"b12x";"a12x",
\ar @{->}"b13x";"a13x",
\ar @{->}"b14x";"a14x",
\POS "a14x",
(142.5,-23.5)*+{\rien}
\endxy}
\end{tabular}
\end{center}
\end{figure}

\begin{figure}
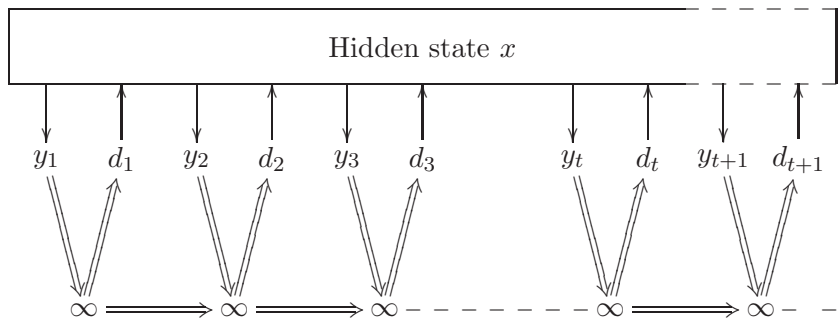

\caption{The optimization process}
\label{ISDA:fig:4}
\begin{center}
\begin{tabular}{c}
\vspace{-25pt}\\
\vspace{-15pt}\rien\hspace{-40pt}\scalebox{1}{\xy
(-12.5,2.5)*+{\rien},
(0,0)="box00",
(90,0)="box10",
(90,-10)="box11",
(110,0)="box10f",
(110,-10)="box11f",
(0,-10)="box01",
\ar @{-}"box00";"box10",
\ar @{--}"box10";"box10f",
\ar @{-}"box10f";"box11f",
\ar @{--}"box11";"box11f",
\ar @{-}"box11";"box01",
\ar @{-}"box00";"box01",\POS "box01",
(55,-5)*+{\mbox{Hidden state $x$}},
(5,-10)="a10y",
(25,-10)="a11y",
(45,-10)="a12y",
(75,-10)="a13y",
(95,-10)="a14y",
(15,-10)="a10x",
(35,-10)="a11x",
(55,-10)="a12x",
(85,-10)="a13x",
(105,-10)="a14x",
(5,-20)*+{y_1}="b10y",
(25,-20)*+{y_2}="b11y",
(45,-20)*+{y_3}="b12y",
(75,-20)*+{y_t}="b13y",
(95,-20)*+{y_{t+1}}="b14y",
(15,-20)*+{d_1}="b10x",
(35,-20)*+{d_2}="b11x",
(55,-20)*+{d_3}="b12x",
(85,-20)*+{d_t}="b13x",
(105,-20)*+{d_{t+1}}="b14x",
\ar @{->}"a10y";"b10y",
\ar @{->}"a11y";"b11y",
\ar @{->}"a12y";"b12y",
\ar @{->}"a13y";"b13y",
\ar @{->}"a14y";"b14y",
\ar @{->}"b10x";"a10x",
\ar @{->}"b11x";"a11x",
\ar @{->}"b12x";"a12x",
\ar @{->}"b13x";"a13x",
\ar @{->}"b14x";"a14x",
\POS "a14x",
(10,-40)*+{\infty}="b00I",
(30,-40)*+{\infty}="b01I",
(50,-40)*+{\infty}="b02I",
(80,-40)*+{\infty}="b03I",
(100,-40)*+{\infty}="b04I",
(110,-40)="b05I",
\ar @2{->}"b00I";"b01I",
\ar @2{->}"b01I";"b02I",
\ar @{--}"b02I";"b03I",
\ar @2{->}"b03I";"b04I",
\ar @2{->}"b10y";"b00I",
\ar @2{->}"b11y";"b01I",
\ar @2{->}"b12y";"b02I",
\ar @2{->}"b13y";"b03I",
\ar @2{->}"b14y";"b04I",
\ar @2{->}"b00I";"b10x",
\ar @2{->}"b01I";"b11x",
\ar @2{->}"b02I";"b12x",
\ar @2{->}"b03I";"b13x",
\ar @2{->}"b04I";"b14x",
\ar @{--}"b04I";"b05I",
\POS "b05I",
(142.5,-50.5)*+{\rien}
\endxy}
\end{tabular}
\end{center}
\end{figure}

\begin{figure}
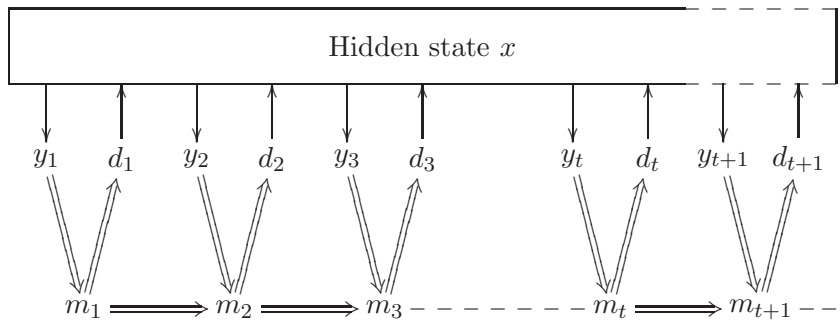

\caption{Finite-memory approximation}
\label{ISDA:fig:5}
\begin{center}
\begin{tabular}{c}
\vspace{-25pt}\\
\vspace{-15pt}\rien\hspace{-40pt}\scalebox{1}{\xy
(-12.5,2.5)*+{\rien},
(0,0)="box00",
(90,0)="box10",
(90,-10)="box11",
(110,0)="box10f",
(110,-10)="box11f",
(0,-10)="box01",
\ar @{-}"box00";"box10",
\ar @{--}"box10";"box10f",
\ar @{-}"box10f";"box11f",
\ar @{--}"box11";"box11f",
\ar @{-}"box11";"box01",
\ar @{-}"box00";"box01",\POS "box01",
(55,-5)*+{\mbox{Hidden state $x$}},
(5,-10)="a10y",
(25,-10)="a11y",
(45,-10)="a12y",
(75,-10)="a13y",
(95,-10)="a14y",
(15,-10)="a10x",
(35,-10)="a11x",
(55,-10)="a12x",
(85,-10)="a13x",
(105,-10)="a14x",
(5,-20)*+{y_1}="b10y",
(25,-20)*+{y_2}="b11y",
(45,-20)*+{y_3}="b12y",
(75,-20)*+{y_t}="b13y",
(95,-20)*+{y_{t+1}}="b14y",
(15,-20)*+{d_1}="b10x",
(35,-20)*+{d_2}="b11x",
(55,-20)*+{d_3}="b12x",
(85,-20)*+{d_t}="b13x",
(105,-20)*+{d_{t+1}}="b14x",
\ar @{->}"a10y";"b10y",
\ar @{->}"a11y";"b11y",
\ar @{->}"a12y";"b12y",
\ar @{->}"a13y";"b13y",
\ar @{->}"a14y";"b14y",
\ar @{->}"b10x";"a10x",
\ar @{->}"b11x";"a11x",
\ar @{->}"b12x";"a12x",
\ar @{->}"b13x";"a13x",
\ar @{->}"b14x";"a14x",
\POS "a14x",
(10,-40)*+{m_1}="b00I",
(30,-40)*+{m_2}="b01I",
(50,-40)*+{m_3}="b02I",
(80,-40)*+{m_t}="b03I",
(100,-40)*+{m_{t+1}}="b04I",
(110,-40)="b05I",
\ar @2{->}"b00I";"b01I",
\ar @2{->}"b01I";"b02I",
\ar @{--}"b02I";"b03I",
\ar @2{->}"b03I";"b04I",
\ar @2{->}"b10y";"b00I",
\ar @2{->}"b11y";"b01I",
\ar @2{->}"b12y";"b02I",
\ar @2{->}"b13y";"b03I",
\ar @2{->}"b14y";"b04I",
\ar @2{->}"b00I";"b10x",
\ar @2{->}"b01I";"b11x",
\ar @2{->}"b02I";"b12x",
\ar @2{->}"b03I";"b13x",
\ar @2{->}"b04I";"b14x",
\ar @{--}"b04I";"b05I",
\POS "b05I",
(142.5,-50.5)*+{\rien}
\endxy}
\end{tabular}
\end{center}
\end{figure}

It is assumed that a subject is acting in a given world with a given purpose or mission.
Thus, the subject interacts with the world and perceives partial informations.
The goal is to optimize the accomplishment of the mission, which is characterized by its reward.
The forthcoming paragraphs are formalizing what is actually a world, what is a mission reward, and how is defined an optimal policy for such a mission.
\paragraph{The world.}
The world is described by an hidden state $x$, which evolves with the time $t$;
in this paper, the time is discretized and increases from step $1$ to step $T$.
More specifically, the variable $x_t$ contains an information which characterizes entirely the world at time $t$.
\emph{In the example of section~\ref{Isda:example}, the hidden state is characterized by the locations of the target and patrols.}
The evolution of the hidden state is given by the vector $x=x_{1:T}=x_1,\dots,x_t,\dots,x_T$.
During the mission, the subject produces decisions $d=d_{1:T}$ which will impact the evolution of the world.
\emph{In example~\ref{Isda:example}, $d$ is the move of the patrols.}
The subject perceives partial observations from the world, denoted $y=y_{1:T}$, which are noisily derived from the hidden state.
\emph{In the example, this observation is an inaccurate estimate of the target location.}
As a conclusion, the world is characterized by a law describing the hidden states and observations conditionnally to the decisions.
This \emph{probabilistic} law is denoted $P$:
\begin{quote}
The hidden state $x_t$ and observation $y_t$ are obtained from the law $P(x_t,y_t|x_{1:t-1},y_{1:t-1},d_{1:t-1})$\,, which are  conditionned by the past hidden states, observations and decisions.
\emph{It is assumed that $d_t$ is generated by the subject after receiving $y_t$\,.}
\end{quote}
In this paper, the law $P$ is quite general, and for example there is no Markovian hypothesis (this hypothesis is required for a dynamic programming approach).
Nevertheless, it is assumed that $P(x_t,y_t|x_{1:t-1},d_{1:t-1})$ may be sampled very quickly.
The law $P(x,y|d)$ is illustrated by figure~\ref{ISDA:fig:1}\,.
In this figure, the out-going arrows are related to the data produced by the world, \emph{i.e.} observations, while incoming arrows are for the data consummed by the world, \emph{i.e.} the decisions.
The variables are put in chronological order from left to right: $y_t$ happens before $d_t$ since decision $d_t$ is produced after observing $y_t$\,.
From now on, $P(x,y|d)$ denotes the law of the world for the completed mission:
$$
P(x,y|d)=\prod_{t=1}^T P(x_t,y_t|x_{1:t-1},y_{1:t-1},d_{1:t-1})\;.
$$
\paragraph{Reward and optimal planning.}
The mission is limited in time and is characterized by a reward.
This reward, denoted $V(d,y,x)$, is a function of the trajectories $d,y,x$\,.
Typically, the function $V$ could be used for computing the time needed for the mission accomplishment.
\emph{The only hypothesis about $V$ is that it is quickly computable.}
In particular, the additivity of the reward\footnote{Additive rewards are of the form $V(d,y,x)=\sum_{t=1}^TV_t(d_t,y_t,x_t)$} with time, a requested hypothesis for many classical methods, is not necessary.
\\[5pt]
The purpose is to construct an optimal decision tree $y\mapsto\bigl(d_t(y_{1:t})|_{t=1}^T\bigr)$, depending on the past observations, in order to maximize the \emph{mean reward}:
\begin{equation}
\label{Isda:Eq:1}
d_\ast\in\arg\max_{d}\sum_{y}\sum_{x}P\left(x,y\left|\bigl(d_t(y_{1:t})|_{t=1}^T\bigr)\right.\right)
V(\bigl(d_t(y_{1:t})|_{t=1}^T\bigr),y,x)\;.
\end{equation}
This optimization process is illustrated by figure~\ref{ISDA:fig:4}.
The double arrows are related to the variables to be optimized.
These arrows describe the information flow between observations and decisions.
\emph{The cells denoted $\infty$ are making decisions and transmitting all the received and generated informations}.
This architecture illustrates that planning with observation is a non-finite memory problem\,: the decision depends on the whole past observations.
Since the optimum for such a problem is generally intractable, it is necessary to search for near-optimal solutions.
The alternative method proposed now relies on the optimal tuning of a probabilitic model of the policies.
\paragraph{Approximating the decision tree.}
In a program like~(\ref{Isda:Eq:1})\,, the variable to be optimized, $d_o$\,, is a \emph{deterministic} object.
In this precise case, $d_o$ is a tree of decision, that is a function which maps to a decision $d_t$ from any sequence of observation $y_{1:t-1}$\,.
But it is more interesting to have a probabilistic viewpoint, when approximating.
Then the problem is equivalent to finding $\pi(d|y)$\,, a probabilistic law of the decisions conditionally to the \emph{past} observations, which maximizes the mean reward:
$$
V(\pi)=\sum_{d}\sum_{y}\sum_{x}\;\;\prod_{t=1}^T\pi(d_t|d_{1:t-1},y_{1:t}) P(x,y|d)\;V(d,y,x)\;.
$$
This new problem is still illustrated by figure~\ref{ISDA:fig:4}, but the double arrows are now describing a Bayesian network structure for the law $\pi$.
By the way, there is not a great difference with the deterministic case for the optimum:
when $d_o$ is unique, the optimal law $\pi_o\in\arg\max_\pi V(\pi)$ is a dirac on $d_o$\,.
However, the probabilistic viewpoint is more suitable to an approximation:
it is simplier to handle probabilistic models than deterministic decision trees, and the optimization is ensured to be continuous; moreover, a natural approximation of $\pi_o$ is obtained by replacing the non-finite memories $\infty$ by finite memories $m$; \emph{c.f.} figure~\ref{ISDA:fig:5}.
Restricting the memory size of the policies is equivalent to approximate the law $\pi$ by a hidden Markov Model.
Then, the approach developped in this paper is quite general and can be split up into two processes:
\begin{itemize}
\item Define a family of parameterized HMMs $\mathcal{H}$\,,
\item Optimize the parameters of the HMM in order to maximize the mean reward:
$$
\mathrm{Find }\quad h_O\in\arg\max_{h\in\mathcal{H}}V(h)\;.
$$
\end{itemize}
As will be seen later, it is easy to tune a HMM optimally by the Cross-Entropy method of Rubinstein\cite{CEbook}.
But first, it is discussed in the next section about the choice of the familly $\mathcal{H}$.
\section{Models}
\begin{figure}
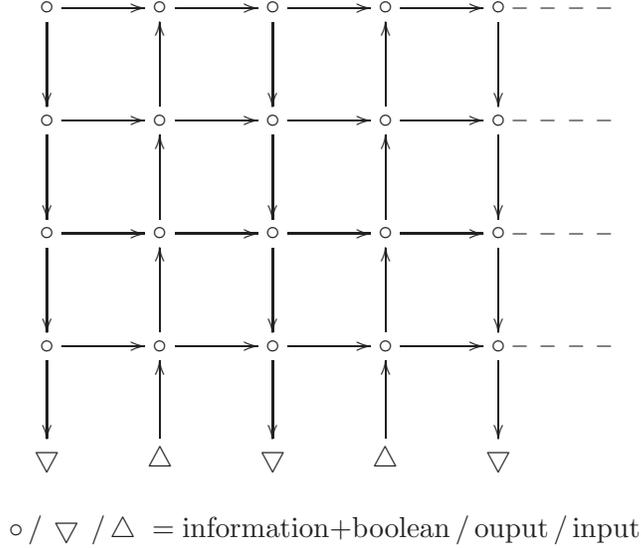

\caption{Model of a controlled Hierarchical HMM}
\label{ISDA:fig:6}
\begin{center}
\begin{tabular}{@{}c@{}}
\vspace{-15pt}\\
\xy
(-.5,.5)*+{},
(0,0)*+{\circ}="b00",
(15,0)*+{\circ}="b01",
(30,0)*+{\circ}="b02",
(45,0)*+{\circ}="b03",
(60,0)*+{\circ}="b04",
(75,0)="b05",
\ar @{->}"b00";"b01",
\ar @{->}"b01";"b02",
\ar @{->}"b02";"b03",
\ar @{->}"b03";"b04",
\ar @{--}"b04";"b05",\POS "b05",
(0,-15)*+{\circ}="b10",
(15,-15)*+{\circ}="b11",
(30,-15)*+{\circ}="b12",
(45,-15)*+{\circ}="b13",
(60,-15)*+{\circ}="b14",
(75,-15)="b15",
\ar @{->}"b00";"b10",
\ar @{<-}"b01";"b11",
\ar @{->}"b02";"b12",
\ar @{<-}"b03";"b13",
\ar @{->}"b04";"b14",
\ar @{->}"b10";"b11",
\ar @{->}"b11";"b12",
\ar @{->}"b12";"b13",
\ar @{->}"b13";"b14",
\ar @{--}"b14";"b15",\POS "b15",
(0,-30)*+{\circ}="b00",
(15,-30)*+{\circ}="b01",
(30,-30)*+{\circ}="b02",
(45,-30)*+{\circ}="b03",
(60,-30)*+{\circ}="b04",
(75,-30)="b05",
\ar @{<-}"b00";"b10",
\ar @{->}"b01";"b11",
\ar @{<-}"b02";"b12",
\ar @{->}"b03";"b13",
\ar @{<-}"b04";"b14",
\ar @{->}"b00";"b01",
\ar @{->}"b01";"b02",
\ar @{->}"b02";"b03",
\ar @{->}"b03";"b04",
\ar @{--}"b04";"b05",\POS "b05",
(0,-45)*+{\circ}="b10",
(15,-45)*+{\circ}="b11",
(30,-45)*+{\circ}="b12",
(45,-45)*+{\circ}="b13",
(60,-45)*+{\circ}="b14",
(75,-45)="b15",
\ar @{->}"b00";"b10",
\ar @{<-}"b01";"b11",
\ar @{->}"b02";"b12",
\ar @{<-}"b03";"b13",
\ar @{->}"b04";"b14",
\ar @{->}"b10";"b11",
\ar @{->}"b11";"b12",
\ar @{->}"b12";"b13",
\ar @{->}"b13";"b14",
\ar @{--}"b14";"b15",\POS "b15",
(0,-60)*+{\bigtriangledown}="b00",
(15,-60)*+{\bigtriangleup}="b01",
(30,-60)*+{\bigtriangledown}="b02",
(45,-60)*+{\bigtriangleup}="b03",
(60,-60)*+{\bigtriangledown}="b04",
\ar @{<-}"b00";"b10",
\ar @{->}"b01";"b11",
\ar @{<-}"b02";"b12",
\ar @{->}"b03";"b13",
\ar @{<-}"b04";"b14",\POS "b14",
(75.5,-60.5)*+{}
\endxy\vspace{10pt}\\
$\circ\,/\,\bigtriangledown\,/\,\bigtriangleup\;=\;$information+boolean\,/\,ouput\,/\,input\vspace{-10pt}
\end{tabular}
\end{center}
\end{figure}

\begin{figure}
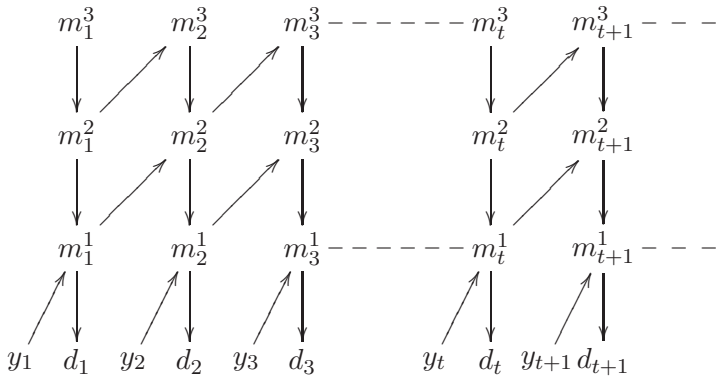

\caption{HHMM model for the planning}
\label{ISDA:fig:7}
\begin{center}
\begin{tabular}{@{}c@{}}
\rien\vspace{-15pt}\\
\xy
(-.5,.5)*+{},
(0,0)*+{m^3_1}="b00",
(15,0)*+{m^3_2}="b01",
(30,0)*+{m^3_3}="b02",
(55,0)*+{m^3_t}="b03",
(70,0)*+{m^3_{t+1}}="b04",
(85,0)="b05",
(0,-15)*+{m^2_1}="b10",
(15,-15)*+{m^2_2}="b11",
(30,-15)*+{m^2_3}="b12",
(55,-15)*+{m^2_t}="b13",
(70,-15)*+{m^2_{t+1}}="b14",
(85,-15)="b15",
\ar @{->}"b00";"b10",
\ar @{->}"b01";"b11",
\ar @{->}"b02";"b12",
\ar @{->}"b03";"b13",
\ar @{->}"b04";"b14",
\ar @{->}"b10";"b01",
\ar @{->}"b11";"b02",
\ar @{->}"b13";"b04",
\ar @{--}"b02";"b03",
\ar @{--}"b04";"b05",\POS "b05",
(0,-30)*+{m^1_1}="b00",
(15,-30)*+{m^1_2}="b01",
(30,-30)*+{m^1_3}="b02",
(55,-30)*+{m^1_t}="b03",
(70,-30)*+{m^1_{t+1}}="b04",
(85,-30)="b05",
\ar @{<-}"b00";"b10",
\ar @{<-}"b01";"b11",
\ar @{<-}"b02";"b12",
\ar @{<-}"b03";"b13",
\ar @{<-}"b04";"b14",
\ar @{->}"b00";"b11",
\ar @{->}"b01";"b12",
\ar @{->}"b03";"b14",
\ar @{--}"b02";"b03",
\ar @{--}"b04";"b05",\POS "b05",
(0,-45)*+{d_1}="b10",
(15,-45)*+{d_2}="b11",
(30,-45)*+{d_3}="b12",
(55,-45)*+{d_t}="b13",
(70,-45)*+{d_{t+1}}="b14",
(85,-45)="b15",
(-7.5,-45)*+{y_1}="b100b",
(7.5,-45)*+{y_2}="b10b",
(22.5,-45)*+{y_3}="b11b",
(47.5,-45)*+{y_t}="b12b",
(62.5,-45)*+{y_{t+1}}="b13b",
\ar @{->}"b00";"b10",
\ar @{->}"b01";"b11",
\ar @{->}"b02";"b12",
\ar @{->}"b03";"b13",
\ar @{->}"b04";"b14",
\ar @{->}"b100b";"b00",
\ar @{->}"b10b";"b01",
\ar @{->}"b11b";"b02",
\ar @{->}"b12b";"b03",
\ar @{->}"b13b";"b04",\POS "b04",
(85.5,-45.5)*+{}
\endxy\vspace{-10pt}
\end{tabular}
\end{center}
\end{figure}

\paragraph{General points.}
The choice of the family of policy models, $\mathcal{H}$, will profoundly impact the efficiency of the approximation.
In particular, the models will be characterized by the memory size and the internal structure of the HMMs (\emph{e.g.} is it hierarchical or not?).
Both characteristics will act upon the convergence, as will be seen in the experiments.
In the most simple case, the HMMs of $\mathcal{H}$ contain no structure and are distinguished by their memory size only.
Example of simple HMM:
\begin{quote}
Let $M$ be indeed a finite set of states, describing the memory capacity of our models.
Then, the memory of the HMM at time $t$ is $m_t\in M$, a variable valued within $M$.
A HMM $h\in\mathcal{H}$ is thus typically defined by:
$$\left\{\begin{array}{@{}l@{}}\displaystyle
h(d|y)=\sum_{m\in M^T}h(d,m|y)\;,
\\\displaystyle
h(d,m|y)=\prod_{t=1}^T \bigl(h_d(d_t|m_t)h_m(m_t|y_t,m_{t-1})\bigr)\;,
\end{array}\right.$$
where the conditionnal law $h_d$ and $h_m$ are time invariant.
\end{quote}
But subsequently will be considered the impact of both the memory and HMM stuctures.
For this purpose a specific family of hierarchical HMM will be introduced and studied.
HHMM are indeed a particular case of HMM, implementing strong intern structures.
\paragraph{Hierarchical HMM.}
Hierarchical models are inspired from biology: to solve a complex problem, factorize it and make decisions in a hierarchical fashion.
Low hierarchies manipulate low level informations and actions, making short-term decisions.
High hierarchies manipulate high level informations and actions (uncertainty is less), making long-term decisions.
Hierarchical HMM are such kind of models.
\emph{A hierarchical hidden Markov model (HHMM) is a HMM which output is either a hierarchical HMM or an actual output.}
A HHMM could also be considered as a hierarchy of \emph{stochastic} processes calling sub-processes.
From this common definition, HHMM are complex structures, which are difficult to formalize and to computerize.
Nevertheless, these models have been introduced and applied for handwriting recognition \cite{IEEEisda:fine}, as well for modelling complex worlds in control applications \cite{IEEEisda:theo}.
A fundamental contribution has been made by \emph{Murphy and Paskin} \cite{IEEEisda:murphy}, which have shown how HHMM could be interpreted as a particular $2-$dimension dynamic Bayesian Network.
Now, Dynamic Bayesian Networks are easily formalized, manipulated and computerized.
DBN could be considered as HMM with complex intern structures.
From the work of \emph{Murphy and Paskin}, it could be shown that a hierarchical HMM (with input and output) could be interpreted by a DBN as described in figure~\ref{ISDA:fig:6}, with discrete or semi-continuous states.
It appears, that there is a up and down flow of the information between the hierarchical levels in addition to the usual temporal flow (the Markovian property).
\emph{It is important to note that boolean informations are necessary for implementing the hierarchy.}
These boolean are needed for controlling the information flows betwenn processes and subprocesses.
The next paragraph introduces the customized model of HHMM, which has been considered in this work.
It is simplification of the general HHMM model, and it allows a more simple implementation.
\paragraph{Implemented model.}
The implemented model familly $\mathcal{H}$ is composed by HHMM with $\Lambda$ hierarchical levels.
Each level $\lambda\in[\![1,\Lambda]\!]$ is associated to a finite memory set $M^\lambda$ (the memory size may change with the hierarchy).
The exchange of information between the levels is characterized by the DBN illustrated in figure~\ref{ISDA:fig:7}.
Notice that each memory cell receives an information from the current upper-level cell and the previous lower-level cell.
As a consequence, the hierarchical and temporal information exchanges are guaranted.
In a more formal way, the HHMM $h\in\mathcal{H}$ are of the form:
$$\left\{\begin{array}{@{}l@{}}\displaystyle
h(d|y)=\sum_{m\in M^{\Lambda T}}h(d,m|y)\;,
\\\displaystyle
h(d,m|y)=\prod_{t=1}^Th^0(d_t|m^1_{t})h^1(m_t^1|y_{t},m^2_{t})
\prod_{\lambda=2}^{\Lambda}h^\lambda(m_{t}^\lambda|m_{t-1}^{\lambda-1},m_{t}^{\lambda+1})\;,
\end{array}\right.$$
where $m^\lambda\in M^\lambda$ is the variable for the memory at level $\lambda$.
It is noteworthy that this model is equivalent to a simple HMM when $\Lambda=2$\,.
And when $\Lambda=1$\,, the law $h$ just maps the immediate observation to decisions, without any memory of the past observations.
\\[5pt]
For any $h\in\mathcal{H}$, define $P[h]$ the complete probabilistic law of the system world/subject:
$$
P[h](d,y,x,m)=P(y,x|d)h(d,m|y)
$$
Then the issue is to find the near-optimal strategy $h_o\in\mathcal{H}$ such that:
$$
h_o\in\arg\max_{h\in\mathcal{H}}\sum_{d,y,x,m}P[h](d,y,x,m) V(d,y,x)\;.
$$
A solution to this problem, by means of the cross-entropy method, is proposed in the next section.
\section{Cross-entropic optimization of $h$}
The reader interested in CE methods should refer to the tutorial \cite{IEEEisda:boer} and the book \cite{CEbook} on the CE method.
CE algorithms were first dedicated to estimating the probability of rare events.
A slight change of the basic algorithm made it also good for optimization.
In their new article\cite{MCO:mello}, Homem-de-Mello and Rubinstein have given some results about the global convergence.
In order to ensure such convergence, some refinements are introduced particularly about the selective rate.
\vspace{5pt}\\
This presentation is restricted to the basic CE method.
The new improvements of the CE algorithm proposed in~\cite{MCO:mello} have not been implemented, but the algorithm has been seen to work properly.
For this reason, this paper does not deal with the choice of the selective rate.
\subsection{General CE algorithm for the optimization}
The Cross Entropy algorithm repeats until convergence the three successive phases:
\begin{enumerate}
\item Generate samples of random data according to a parameterized random mechanism,
\item Select the best samples according to a reward criterion,
\item Update the parameters of the random mechanism, on the basis of the selected samples.
\end{enumerate}
In the particular case of CE, the update in phase 3 is obtained by minimizing the Kullback-Leibler distance, or cross entropy, between the updated random mechanism and the selected samples.
The next paragraphs describe on a theoretical example how such method can be used in an optimization problem.
\paragraph{Formalism.}
Let be given a function $x\mapsto f(x)$; this function is easily computable.
The value $f(x)$ has to be maximized, by optimizing the choice of $x\in X$.
The function $f$ will be the reward criterion.
\vspace{5pt}\\
Now let be given a family of probabilistic laws, $P_\sigma|_{\sigma\in\Sigma}$\,, applying on the variable $x$.
The family $P$ is the parameterized random mechanism.
The variable $x$ is the random data.
\vspace{5pt}\\
Let $\rho\in\,]0,1[$ be a selective rate.
The CE algorithm for $(x,f,P)$ follows the synopsis\,:
\begin{enumerate}
\item Initialize $\sigma\in\Sigma$\,,
\item \label{XX:step2}Generate $N$ samples $x_n$ according to $P_\sigma$\,,
\item Select the $\rho N$ best samples according to the reward criterion $f$\,,
\item Update $\sigma$ as a minimizer of the cross-entropy with the selected samples:
$$\sigma\in\arg\max_{\sigma\in\Sigma}\sum_{n~{\rm selected}}\ln P_\sigma(x_n)\;,$$
\item Repeat from step \ref{XX:step2} until convergence.
\end{enumerate}
\emph{This algorithm requires $f$ to be easily computable and the sampling of $P_\sigma$ to be fast.}
\paragraph{Interpretation.}
The CE algorithm tightens the law $P_\sigma$ around the maximizer of $f$.
Then, when the probabilistic family $P$ is well suited to the maximization of $f$\,, it becomes equivalent to find a maximizer for $f$ or to optimize the parameter $\sigma$ by means of the CE algorithm.
The problem is to find a good family\dots
Another issue is the criterion for deciding the convergence.
Some answers are given in \cite{MCO:mello}.
Now, it is outside the scope of this paper to investigate these questions precisely.
Our criterion was to stop after a given threshold of successive \emph{unsuccessful tries} and this very simple method have worked fine on our problem.
\subsection{Application}\label{MCO:Subsect:MainAlgo}
Optimizing $h\in\mathcal{H}$ means tuning the parameter $h$ in order to tighten the probability $P[h]$ around the optimal values for $V$\,.
This is exactly solved by the \emph{Cross-Entropy} optimization method.
However, it is required that the reward function $V$ is easily computable.
Typically, the definition of $V$ may be recursive, \emph{e.g.}\,:
$$
V(d,y,x)=V_T\ ; \quad V_{t}=v_{t}(d_t,y_t,x_t,V_{t-1})\quad\mbox{and}\quad V_0=0\;.
$$
Let the \emph{selective rate} $\rho$ be a positive number such that $\rho<1$\,.
The cross-entropy method for optimizing $h$ follows the synopsis\,:
\begin{enumerate}
\item Initialize $h$\,. For example a flat $h$,
\item \label{Algo:0:1:FirstStep}Build $N$ samples $\theta^n=(d^n,y^n,x^n,m^n)$ according to the law $P[h]$,
\item Choose the $\rho N$ best samples $\theta^n$ according to the reward $V(d^n,y^n,x^n)$\,.
Denote $S$ the set of the selected samples,
\item Update $h$ as the minimizer of the cross-entropy with the selected samples:
\begin{equation}
\label{StratApp:1:1}
h\in\arg\max_{h\in\mathcal{H}}\sum_{n\in S}\ln P[h](\theta^n)\;,
\end{equation}
\item Reiterate from step~\ref{Algo:0:1:FirstStep} until convergence.
\end{enumerate}
For our HHMM model, the maximization~(\ref{StratApp:1:1}) is solved by:
$$
h^0(A|B)=\frac{
\mathrm{card}\Bigl\{n\in S\,,t\,/\,A=d_t^{n}\;~\mbox{and}~
B=m_t^{1,n}\Bigr\}
}{
\mathrm{card}\Bigl\{n\in S\,,t\,/\,B=m_t^{1,n}\Bigr\}
}\;,
$$
$$
h^1(A|B,C)=\frac{
\mathrm{card}\Bigl\{n\in S\,,t\,/\,A=m_t^{1,n}\,,
B=y_{t}^{n}\;~\mbox{and}~C=m_t^{2,n}\Bigr\}
}{
\mathrm{card}\Bigl\{n\in S\,,t\,/\,B=y_{t}^{n}\;~\mbox{and}~C=m_t^{2,n}\Bigr\}
}\;.
$$
and for $2\le\lambda\le\Lambda$\,,:
$$
h^\lambda(A|B,C)=\frac{
\mathrm{card}\Bigl\{n\in S\,,t\,/\,A=m_t^{\lambda,n}\,,
B=m_{t-1}^{\lambda-1,n}\;~\mbox{and}~C=m_t^{\lambda+1,n}\Bigr\}
}{
\mathrm{card}\Bigl\{n\in S\,,t\,/\,B=m_{t-1}^{\lambda-1,n}\;~\mbox{and}~C=m_t^{\lambda+1,n}\Bigr\}
}\;.
$$
The next section presents an example of implementation of the algorithm described in section~\ref{MCO:Subsect:MainAlgo}.
\section{Implementation}
\label{Isda:example}
The algorithm has been applied to a simulated target detection problem.
\subsection{Problem setting}
A target $R$ is moving in a lattice of $20\times 20$ cells, \emph{ie.} $[\![\,0,19\,]\!]^2$.
$R$ is tracked by two mobiles, $B$ and $C$, controlled by the subject.
The coordinate of $R$, $B$ and $C$ at time $t$ are denoted $(i_R^t,j_R^t),$ $(i_B^t,j_B^t)$ and $(i_C^t,j_C^t)$.
$B$ and $C$ have a very limited information about the target position, and are maneuvering much slower:
\begin{itemize}
\item A move for $B$ (respectively $C$) is either: \emph{turn left}, \emph{turn right}, \emph{go forward}, \emph{no move}.
Consequently, there are $4\times4=16$ possible actions for the subject.
These moves cannot be combined in a single turn. No diagonal forward: a mobile is either directed up, right, down or left,
\item The mobiles are initially positioned in the down corners, \emph{ie.} $i^1_B=0,$ $j^1_B=19$ and $i^1_C=19,$ $j^1_C=19$.
The mobile are initially directed \emph{downward},
\item $B$ (respectively $C$) observes whether the target relative position is forward or not.
More precisely:
\begin{itemize}
\item when $B$ is directed upward, it knows whether $j_R<j_B$ or not,
\item when $B$ is directed right, it knows whether $i_R>i_B$ or not,
\item when $B$ is directed downward, it knows whether $j_R>j_B$ or not,
\item when $B$ is directed left, it knows whether $i_R<i_B$ or not,
\end{itemize}
\item $B$ (respectively $C$) knows whether its distance with the target is less than $3$, \emph{ie.} $d_\infty(B,R)<3$, or not.
The distance $d_\infty$ is defined by:
$$
d_\infty(B,R)=\max\{|i_B-i_R|\,,|j_B-j_R|\}\;.
$$
\end{itemize}
At last, there are $2^4=16$ possible observations for the subject.
\\[5pt]
Several test cases have been considered.
In case 1, the target $R$ does not move.
In any other case, the target $R$ chooses stochastically its next position in its neighborhood.
Any move is possible (up/down, left/right, diagonals, no move).
The probability to choose a new position is proportional to the sum of the squared distance from the mobiles:
$$
\left\{
\begin{array}{@{}l@{}}\displaystyle
P(R^{t+1}|R^t)=0\mbox{ if }|i_{R}^{t+1}-i_{R}^{t}|>1\mbox{ or }|j_{R}^{t+1}-j_{R}^{t}|>1\;, \vspace{5pt}
\\\displaystyle
P(R_{t+1}|R_t)\propto(i_{R}^{t+1}-i_{B}^{t})^2+(j_{R}^{t+1}-j_{B}^{t})^2
\\\displaystyle\hspace{170pt}
+(i_{R}^{t+1}-i_{C}^{t})^2+(j_{R}^{t+1}-j_{C}^{t})^2\quad\mbox{else}\;.
\end{array}\right.
$$ 
This definition was intended to favorize escape moves: more great is a distance, more probable is the move.
But in such summation, a short distance will be neglected compared to a long distance.
It is implied that a distant mobile will hide a nearby mobile.
This ``deluding'' property will induce actually two different kinds of strategy, whithin the learned machines.
\\[5pt]
The objective of the subject is to maintain the target sufficiently closed to at least one mobile (in this example, the distance between the target and a mobile is required to be not more than $3$).
More precisely, the reward function, $V$, is just counting the number of such ``encounter'':
$$
V_0=0\ ;\quad V_{t}=V_{t-1}+1\mbox{\ \ if }d_\infty(B^t,R^t)\le3\mbox{ or }d_\infty(C^t,R^t)\le3\ ;\quad
V_{t}=V_{t-1}\mbox{\ \ else.}
$$
The total number of turns is $T=100$.
\subsection{Results}
\paragraph{Generality.}
Like many stochastic algorithms, this algorithm needs some time for convergence.
For the considered example, about two hours were needed for convergence (on a 2GHz PC); the selective rate was $\rho=0.5$.
This speed depends on the size of the HHMM model and on the convergence criterion.
A weak and a strong criterion are used for deciding the convergence.
Within the weak criterion, the algorithm is terminated after $100$ successive unsuccessful tries.
Within the strong criterion, the algorithm is terminated after $500$ successive unsuccessful tries.
Of course, the strong criterion computes a (slightly) better optimum than the weak criterion, but it needs time.
Because of the many tested examples, the weak criterion has been the most used in particular for the big models.
For the same HHMM model, the computed optimal values do not depend on the algorithmic instance (small variations result however from the stochastic nature of the algorithm).
\\[5pt]
In the sequel, mean rewards are rounded to the nearest integer, or are expressed as a percentage of the optimum.
Thus, the presentation is made clearer.
And owing to the small variations of this stochastic algorithm, more precision turns out to be irrelevant.
\paragraph{Case 1: $R$ does not move.}
This example has been considered in order to test the algorithm.
The position of the target is fixed in the center of the square space, \emph{ie.} $i^1_R=j^1_R=10$.
It is recalled that the mobiles are initially directed downward.
Then, the optimal strategy is known and its value is $85$\,: the time needed to reach the target is $15$\,, and no further move is needed.
The learned $h_o$ approximates the reward $84$\,.
The convergence is good.
\paragraph{Case 2: $R$ is moving but the observation $y$ is hidden.}
Initially, $R$ is located within the $20\times10$ upper cells of the lattice (\emph{ie.} $[\![\,0,19\,]\!]\times[\![\,0,9\,]\!]$), accordingly to a uniform probabilistic law.
The computed optimal means reward is about $32$.
In this case, the mobiles tend to move towards the upper corners.
\paragraph{Case 3: $R$ is moving and $y$ is observed.}
Again, $R^1$ is located uniformly within the $20\times10$ upper cells of the lattice.
The computed optimal means reward is about $69$.
This reward has been obtained from a large HHMM model ($\Lambda=2$ with $256$ states per level, \emph{ie.} $\mathrm{card}(M^\lambda)=256$) and with the strong criterion.
However, somewhat smaller models should work as well.
\\[3pt]
Specific computations are now presented, depending on the number of levels $\Lambda$ and the number of states per levels.
For each case, the weak criterion has been used.
The rewards are now expressed as percentage.
\\[5pt]
\emph{Subcase $\Lambda=1$\,.}
For such model, the action $d_t$ is constructed only from the immediate last observation $y_{t}$.
The model does not keep any memory of the past observations.
Then, only $16$ states are sufficient to describe the hidden variable $m_t^1$\,, \emph{ie.} $\mathrm{card}(M^1)=16$.
The resulting reward is $78\%$ of the optimum.
\\[5pt]
\emph{Subcase $\Lambda=2$\,.}
This model is equivalent to a HMM and it is assumed that $\mathrm{card}(M^1)=\mathrm{card}(M^2)$.
The following table gives the computed reward for several choices of the memory size:
$$
\begin{array}{|l||l|l|l|l|}
\hline \mathrm{card}(M^\lambda)&16&32&64&256
\\
\hline\mathrm{Reward}&94\%&96\%&97\%&97\%
\\\hline
\end{array}
$$
%
It is noteworthy that the memory of the past observations allows better strategies than the only last observation (case $\Lambda=1$)\,.
Indeed, the reward jumps from $78\%$ up to $97\%$.
\\[5pt]
\emph{Subcases $\Lambda>2$.}
A comparison of graduated hierarchic models, $1\le\Lambda\le4$, has been made.
The first level contained $16$ possible states, and the higger levels were restricted to $2$ states:
$$
\begin{array}{|l||l|l|l|l|}
\hline\mathrm{hierarchic\ grade}&\Lambda=1&\Lambda=2&\Lambda=3&\Lambda=4
\\
\hline\mathrm{card}(M^\lambda)|_{\lambda=1}^\Lambda&16&16,2&16,2,2&16,2,2,2
\\\hline
\end{array}
$$
The test has been accomplished according to the weak criterion:
$$
\begin{array}{|l||l|l|l|l|}
\hline\mathrm{hierarchic\ grade}&\Lambda=1&\Lambda=2&\Lambda=3&\Lambda=4
\\
\hline\mathrm{Reward\ (weak)}&78\%&85\%&\bf 81\%&94\%
\\\hline
\end{array}
$$
and the strong criterion:
$$
\begin{array}{|l||l|l|l|l|}
\hline\mathrm{hierarchic\ grade}&\Lambda=1&\Lambda=2&\Lambda=3&\Lambda=4
\\
\hline\mathrm{Reward\ (strong)}&80\%&88\%&93\%&96\%
\\\hline
\end{array}
$$
It seems that a high hierarchic grade (\emph{i.e.} more structure) makes the convergence difficult.
This is particularly the case here for the grade $\Lambda=3$\,, which failed under the weak criterion at only $81\%$.
However, the algorithm works again when improving the convergence criterion.
\\[5pt]
It is interesting to make a comparison with the subcase $\Lambda=2$ where $\mathrm{card}(M^1)=\mathrm{card}(M^2)=16$.
Under the weak criterion, the result for this HHMM was $94\%$ as for the grade $\Lambda=4$.
However, the dimension of the law is quite different for the two models:
\begin{itemize}
\item $15\times 16+15\times 16\times 16+15\times 16=4320$ for the $2$-level HHMM,
\item $15\times 16+15\times 16\times 2+1\times 16\times 2+1\times 2\times 2+1\times 2=758$ for the $4$-level HHMM.
\end{itemize}
This dimension is a rough characterization of the complexity of the model.
It seems clear on these examples that the highly hierarchized models are more efficient than the weakly hierarchized models.
And the problem considered here is quite simple.
On complex problems, hierarchical models may be pre-eminent.
\paragraph{Global behavior.}
\rien\\\emph{The algorithm.}
The convergence speed is low at the beginning.
After this initial stage, it improves greatly until it reaches a new ``waiting'' stage.
This alternation of low speed and great speed stages have been noticed several times.
\\[5pt]
\emph{The near optimal policy.}
It is now discussed about the behaviour of the best found policy.
This policy has reach the mean reward $69$.
The mobiles strategy results in a tracking of the target.
The figure~\ref{ISDA:fig:8} illustrates a short sequence of escape/tracking of the target.
It has been noticed two quite distinct behaviours, among the many runs of the policy:
\begin{itemize}
\item \emph{The two mobiles may both cooperate on tracking the target,}
\item \emph{When the target is near a border, one mobile may stay along the opposite border while the other mobile may perform the tracking.}
This strategy seems strange at first sight.
But it is recalled that the moving rule of the target tends to neglect a nearby mobile compared to a distant mobile.
In this strategy, the first mobile is just annihilating the ability of the target to escape from the tracking of the second mobile.
\end{itemize}
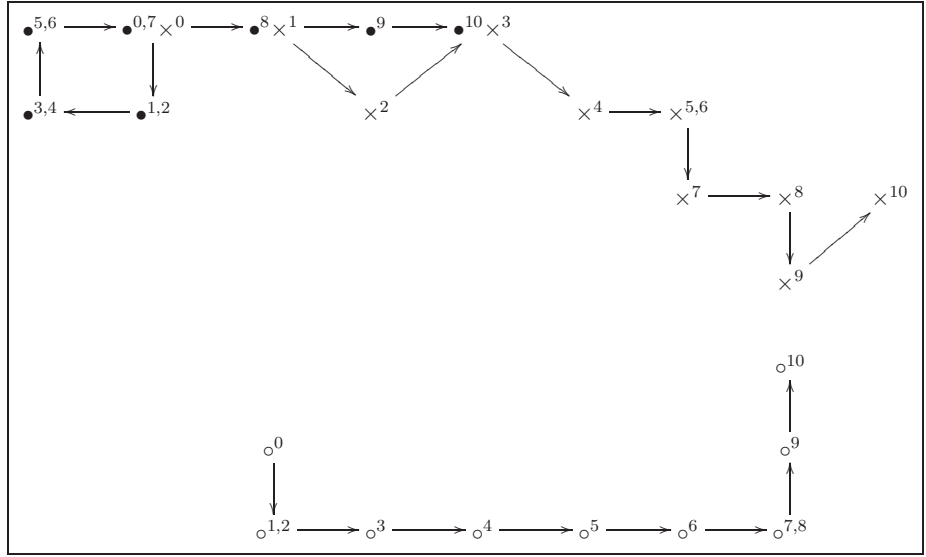
\begin{figure}[ht!]
\caption{Near-optimal control sequence}
\label{ISDA:fig:8}
\begin{center}
\begin{tabular}{l}\vspace{-15pt}\\
\fbox{\scalebox{.8}{\xymatrix{
\bullet^{5,6}\ar[r]&\bullet^{0,7}\ar[d]\times^0\ar[r]&\bullet^8\ar[r]\times^1\ar[dr]&\bullet^9\ar[r]&\bullet^{10}\times^3\ar[dr]&&&&\\
\bullet^{3,4}\ar[u]&\bullet^{1,2}\ar[l]&&\times^2\ar[ur]&&\times^4\ar[r]&\times^{5,6}\ar[d]&&\\
&&&&&&\times^7\ar[r]&\times^8\ar[d]&\times^{10}\\
&&&&&&&\times^9\ar[ru]&\\
&&&&&&&\circ^{10}&\\
&&\circ^0\ar[d]&&&&&\circ^9\ar[u]&\\
&&\circ^{1,2}\ar[r]&\circ^{3}\ar[r]&\circ^{4}\ar[r]&\circ^{5}\ar[r]&\circ^{6}\ar[r]&\circ^{7,8}\ar[u]&
}}}\\\\
$\times=\;$target\qquad
$\bullet=\;$observer 1\qquad
$\circ=\;$observer 2
\\
Relative times are put in supscript\vspace{-10pt}
\end{tabular}
\end{center}
\end{figure}
\subsection{Comparison with the Q-learning}
The Q-learning is a reinforcement learning method, which is based on the computation of a table evaluating the decision conditionnally to the \emph{known information}.
The known information is typically the state of the world if it is known, or partial states and observations.
Since the known information increases exponentially with the observation range, the test will only implement a Q-learning based on the immediate past observation.
Now, let us recall some theoretical grounds about the Q-learning.
\paragraph{Theory.}
A founding reference about reinforcement learning is the well known book of Sutton and Barto \cite{sutton}, which is available on internet.
This paragraph will not enter deeply into the subject, and is limited to a simple description of the Q-learning.
Moreover, we will make the hypothesis of infinite horizon (that is $T=\infty$) with a weak discounting of the reward $\gamma=0.99$, so as to implement the algorithm in its most classical form.
Tests however have also been made with a finite horizon but have not achieved a good convergence for the considered algorithm.
\\[5pt]
The learning relies on the following hypotheses:
\begin{itemize}
\item At each step $t$, the subject has a (partial) knowledge $s$ of the state of the world, and chooses an action $a$,
\item Let $V_{t+1}$ be the cumulated reward from step $t+1$ to step $\infty$.
Assume a state $s_t$ and action $a_t$ at step $t$.
Then $V_t=R(s_t,a_t)+\gamma V_{t+1}$\,, \emph{i.e.} an instantaneous reward $R$ is obtained and cumulated to the discounted future reward.
\end{itemize}
The question is: \emph{being given a current state $s$, what is the best action $a$ to be done?}
The answer is simple, if we are able to predict the future and evaluate the expected cumulated reward $Q(s,a)$ for any $a$: the best action is $a_o\in\arg\max_a Q(s,a)$\,.
The following algorithm could be used for learning the table $Q$ (taken from \cite{sutton})\,:
\begin{itemize}
\item Initialize $Q(s,a)$ arbitrary
\item (Repeat for each episode:\ \ \emph{[finite-horizon case]})
\begin{itemize}
\item Initialize $s$
\item Repeat for each step (of the episode):
\begin{itemize}
\item With probability $1-\epsilon$ choose $a\in\arg\max_aQ(s,a)$\,; otherwise chose $a$ randomly
\item Take action $a$, receive reward $R(s,a)$ and observe the new state $s'$
\item Set $Q(s,a):=Q(s,a)+\alpha\bigl(R(s,a)+\gamma\max_{a'}Q(s',a')-Q(s,a)\bigr)$
\item Set $s:=s'$
\end{itemize}
\item (until $s$ is terminal)
\end{itemize}
\end{itemize}
where $\alpha$ controls the convergence speed and $\epsilon$ the innovation.
\\[5pt]
In our implementation, $s=(y_t,i_B^t,j_B^t,i_C^t,j_C^t,\mbox{directions})$, $a=d_t$, $\alpha=0.1$, $\epsilon=1/\ln t$ and the instantaneous reward $R$ is complient with the experiment definition of previous section.
Since $s$ contains the last observation plus the known part of the world state,
this experiment should be equivalent to [case 3/subcase $\Lambda=1$] considered previously. 
The computer memory needed to store the table $Q$ was approximately $2$ giga-byte: we are around the limits of the computer.
In particular, it is rather uneasy to involve a greater observation range without some approximations.
\paragraph{Results.}
The algorithm has been stoped after $10^{11}$ iterations, but $10^{10}$ seemed sufficient.
It took several hours, but the algorithm has not been optimized.
In order to make the comparison possible with our method, the Q-strategies has been evaluated by a non-discounted cumulation of the reward on 100-step-wide windows.
Moreover, these evaluations have been made:
\begin{itemize}
\item from the initial stage of the simulation, so as to conform to previous section,
\item after many cycles, so as to simulate an infinite horizon.
\end{itemize}
The following table makes a comparison between the Q-strategies and the model based strategies with $\Lambda=1$.
$$
\begin{array}{|l||l|l|l|}
\hline &\mbox{worse}&\mbox{mean}&\mbox{best}
\\
\hline\mbox{Q-policy/stage $0$}&0\%&40\%&112\%
\\
\hline\mbox{Q-policy/$\infty$-horizon}&0\%&51\%&145\%
\\
\hline\mbox{Model based }\Lambda=1&44\%&78\%&105\%
\\\hline
\end{array}
$$
It is first noticed that the policy obtained by the Q-learning is less regulated than the model based policy.
Moreover, although it may be quite good to track a target when the encounter has been inited (best is $145\%$), it is rather bad at initing the encounter (mean for initial stage is $40\%$) or when the tracking is lost (worst is $0\%$).
At last, the mean evaluation at infinite horizon is $51\%$, which is even smaller than the model-based policy working from the initial stage.
\\[5pt]
On this example, and for this simple Q-learning implemention, the comparison is favorable to the model-based policy.
Moreover, model-based policies are able to manage more observation range.
Now, this planning example has been constructed so as to make difficult the management of the state variables (the dimension is huge) and observations (the observations are poor and have to be combined).
For such a problem, a more dedicated RL-method should be chosen.
\section{Conclusion}
In this paper, we proposed a general method for approximating the optimal planning in a partially observable world.
Hierarchical HMM families have been used for approximating the optimal decision tree, and the approximation has been optimized by means of the Cross-Entropy method.
\\[5pt]
At this time, the method has been applied to a strictly discrete-state problem and has been seen to work properly.
This algorithm has been compared favorably with a Q-learning implementation of the considered problem: it is able to manage more observation range, and the optimized policy is more regulated.
An interesting point is that the optimized policy has discovered two quite different global strategies and is able to choose between them: make the mobiles both cooperate on tracking or require one mobile for deluding the target.
\\[5pt]
The results are promising.
However, the observation and action spaces are limited to a few number of states.
And what happens if the hidden space becomes much more intricated?
There are several possible answers to such difficulties:
\\[5pt]
First, the cross-entropic principle could be applied for optimizing continuous laws.
It is thus certainly possible to consider semi-continuous models, which will be more realistic for a planning policy.
Secondly, many refinements are foreseeable about the structure of the models.
Hierarchic models for observation, decision and memory should be improved in order to locally factorize intricated problems.
This research is just preliminary and future works should investigate these questions.
%

\end{document}